\documentclass[12pt]{article}
\usepackage{amsmath}
\usepackage{amssymb}
\usepackage{amsthm}
\usepackage{url}
\usepackage[numbers]{natbib}
\usepackage{pgf}
\usepackage{tikz}

\usetikzlibrary{arrows,automata}
\usetikzlibrary[positioning]
\tikzstyle{state}=[circle,minimum size=1.5cm,draw=black,scale=0.6,node distance=4cm]

\theoremstyle{plain}
\newtheorem{theorem}{Theorem}
\newtheorem*{lemma-edgecon}{Lemma 1}
\newtheorem*{lemma-partition}{Lemma 2}
\newtheorem*{lemma-subset}{Lemma 3}
\newtheorem*{theorem-strongn}{Theorem 4}
\newtheorem{lemma}{Lemma}

\newtheorem{conjecture}{Conjecture}
\theoremstyle{definition}
\newtheorem{definition}{Definition}
\newtheorem{corollary}{Corollary}

\numberwithin{equation}{section}

\newcommand{\stirling}[2]{\genfrac{\{}{\}}{0pt}{}{#1}{#2}}
\newcommand{\wordgraphfam}[2]{\ensuremath{\mathfrak{W}(#1,#2)}}

\title{The strong--connectivity of word--representable digraphs}

\author{Edward J. L. Bell\\
\small Department of Maths \& Stats\\[-0.8ex]
\small Lancaster University, UK \\
\small \texttt{e.bell@comp.lancs.ac.uk}\\
\and
Paul Rayson\\
\small School of Computing and Communications\\[-0.8ex]
\small Lancaster University, UK\\[-0.8ex]
\small \texttt{p.rayson@lancs.ac.uk}
\and
Damon Berridge \\
\small Department of Maths \& Stats\\[-0.8ex]
\small Lancaster University \\[-0.8ex]
\small \texttt{d.berridge@lancs.ac.uk}
}

\begin{document}
\maketitle

\begin{abstract}
A word--graph $G_\omega$ is a digraph represented by a word $\omega$ such that the vertex--set $V(G_\omega)$ is the alphabet of $\omega$ and the edge--set $E(G_\omega)$ is determined by non--identical adjacent letter pairs in $\omega$. In this paper we study the strong--connectivity of word--graphs. Our main result is that the number of strongly connected word--graphs represented by $\ell-$words of over an $n-$alphabet can be expressed via the recurrence relation $T(\ell,n)$ on the Stirling numbers of the second kind using a link between word partitions and digraph connectivity.
$$T(\ell, n) = \stirling{\ell-1}{n} + \sum_{j=0}^{\ell-2} \sum_{m=0}^{n-2} \stirling{j}{m} T(\ell-j-1, n-m) (n-m-1)$$
\end{abstract}

\section{Introduction}
A digraph $D$ is representable if $D$ can be mapped to a word in some manner.  In this paper we study a particular form of edge--set encoding for word--representable digraphs that has close links with discrete stochastic processes.
\begin{definition}
\label{def:edges}
A word--graph $G_\omega$ is a simple digraph associated with an $\ell-$word $\omega$ over an $n-$alphabet for $\ell > n > 0$ such that the vertex set of $G_\omega$, $V(G_\omega)$, is the alphabet of $\omega$ and the directed edge--set of $G_\omega$,  $E(G_\omega)$, is encoded from $\omega$ using the set of non--identical adjacent letter pairs:

\begin{align*}
V(G_\omega) &= A(\omega) \\
E(G_\omega) &= E(\omega) = \{ (\omega_1, \omega_2), (\omega_2, \omega_3), \ldots, (\omega_{\ell-1}, \omega_\ell)\} 
\end{align*}
such that for $\{\omega_i,\omega_{i+1}\} \in E(G_\omega)$ it holds that $\omega_i \neq \omega_{i+1}$
\end{definition}

The edge--set encoding in Definition~\ref{def:edges} is of particular note due to its close links with a range of stochastic word--based processes such as language, music, DNA and protein sequences. This form of edge--set encoding is also closely related to the concept of digraph path decomposition\cite{path74}. Other forms of edge--set encoding have been used in the literature for both word--representable graphs\cite{representable, word-tern} and word--representable digraphs\cite{overlap}. These existing encodings are used for very specific applications (such as semigroup or DNA sequence analysis) whilst word--graphs have a wide range of general applications.

The focus of this paper is on determining the criteria required for a word to represent a strongly--connected digraph and how many words, out of all $\ell-$words over an $n-$alphabet, represent strongly--connected digraphs. Strong--connectivity is an vital concept for word--graph theory because every strongly--connected word--graph, or component, can be represented by a single word\cite{wordgraph}. 

\begin{definition}
\label{def:strong}
A digraph $D$ is strongly--connected if there exists a path from every $v_i \in V(D)$ to every $v_j \in V(D)$.
\end{definition}

This paper is organised into four sections. The remainder of the first section details previous results and the notation used in this paper. The second section describes our results and the third section contains proofs of these results. Future work is detailed in the final section.

\subsection{Previous results}
\label{subsec:previous}
We say a digraph $D$ is representable, \emph{i.e.} is a word--graph, if there exists a word $\omega$ that encodes $D$ in the manner specified by Definition~\ref{def:edges}. Theorem~\ref{th:path} provides the criterion to be satisfied for a digraph to be representable although there are other ways to characterise representability such as forbidden substructures or condensation forms\citep{wordgraph}.

\begin{theorem}
\label{th:path}
A digraph $D$ is a word--graph (representable) if there exists a directed path $P$ from $v_i \in V(D)$ to $v_j \in V(D)$ such that $P$ transverses each edge of $D$ at least once. The sequence of vertices visited by $P$ is a representational word of $D$.
\end{theorem}

\begin{corollary}
\label{col:weak}
Representable digraphs are weakly--connected.
\end{corollary}

\begin{corollary}
\label{col:strong}
Strongly--connected digraphs are representable.
\end{corollary}

In the same manner that digraphs are characterised by vertex and edge set cardinality, word--graphs are characterised by word length and alphabet cardinality. The number of word--graphs represented by words of length $\ell$ over an alphabet of cardinality $n$ can be expressed using Stirling numbers of the second kind\cite{word-graph-grammar}.

\begin{theorem}
\label{th:card}
The set of $\ell-$words over an $n-$alphabet is of cardinality $n! \stirling{\ell}{n}$.
\end{theorem}

The term $\stirling{\ell}{n}$ is a Stirling number of the second kind, or the number of ways of partitioning a set of cardinality $\ell$ into $n$ subsets. If $p$ is an $n-$partition of $\ell$ then the word corresponding to $p$ would be: $\omega = p_1p_2\ldots p_n$, for some ordering (\emph{i.e.} lexicographical) of the partition members. The $n!$ term in Theorem~\ref{th:card} accounts for the permutation of the partitions to reflect the order in which they can occur in a word. We refer to the set of $\ell-$words over an $n-$alphabet as a word--graph family $\wordgraphfam{\ell}{n}$. By Corollary~\ref{col:strong}, every strongly--connected digraph is represented by a word in some word--graph family.

Finally we define Menger's Theorem\cite{menger} which is used in the proofs in Section~\ref{sec:proofs}.

\begin{theorem}
The edge--connectivity of a digraph is equal to the minimum number of edge--disjoint paths between any vertex pair.
\end{theorem}

\subsection{Notation and terminology}
\label{sec:not}
An $\ell-$word $\omega$ over an $n-$alphabet is an ordered sequence of symbols of length $\ell = |\omega|$: $\omega = \omega_1\ldots \omega_\ell$ such that $\omega_i \in A(\omega)$ for $|A(\omega)|=n$ where $A(\omega)$ is the alphabet of $\omega$. A circular word is one for which $\omega_1 = \omega_\ell$ and a trivial word is one for which $|A(\omega)|=1$.

A factor of a word $\omega$ is any $\omega_j\ldots\omega_k$ contiguous subword. A partition, or a factorisation, $p$ of a word $\omega$ is a $k-$tuple for $1\le k \le \ell$ of factors such that $\omega = p_1p_2\ldots p_k$. A disjoint partition of $\omega$ is a partition such that, for each $p_i \neq p_j$ pair of partitions, the alphabets of $p_i$ and $p_j$ are disjoint: $A(p_i) \cap A(p_j) = \varnothing$. A minimal disjoint partition of $\omega$ is a partition of $\omega$ with the smallest cardinality such that the partitions remain disjoint. 

The edge--connectivity of a word--graph is $\lambda(G_\omega)$ where $\lambda(G_\omega)$ is the cardinality of the edge--cut set $S \subset E(G_\omega)$ such that $S$ is the smallest set of edges that disconnects $G_\omega$ when removed.

\section{Results}
\label{sec:results}
The initial results in this paper take the form of three lemmas concerning the structure of words that represent strongly--connected digraphs. The three lemmas build upon one another. The final lemma is used to prove a theorem on the number of words that represent strongly--connected word--graphs in a word--graph family. The three lemmas are based upon the link between word--graph strong--connectivity and edge--connectivity.

\begin{lemma}
\label{lem:edgecon}
Iff $\lambda(G_\omega) > 1$ then $G_\omega$ is strongly--connected.
\end{lemma}

A bridge is an unidirectional edge in $G_\omega$ that when removed disconnects $G_\omega$. When such a bridge exists in $G_\omega$ we say $G_\omega$ is 1--edge-connected or that $\lambda(G_\omega)=1$. Lemma~\ref{lem:edgecon} proves that word--graphs with bridges are the only word--graph which are not strongly--connected. This edge--connectivity property holds due to the way in which edges are encoded from the representational word (Definition~\ref{def:edges}). Using Lemma~\ref{lem:edgecon} we prove Lemma~\ref{lem:partition} which characterises the structural properties of words that represent strongly--connected digraphs.

\begin{lemma}
\label{lem:partition}
Iff $\lambda(G_\omega) > 1$ then $\omega$ cannot be partitioned into disjoint factors.
\end{lemma}

See Subsection~\ref{sec:not} for a definition of \emph{disjoint factors}. Lemma~\ref{lem:partition} is restated by Lemma~\ref{lem:subset} in a manner that facilitates the enumeration of strongly--connected word--graphs.

\begin{lemma}
\label{lem:subset}
The number of $\ell-$words over an $n-$alphabet, with a unique factorisation, that represent word--graphs with $\lambda(G_\omega) > 1$ is the number of $n-$partitions of an $\ell-$set such that no proper subset of parts has a union equal to the subset $\{1,2,\ldots,j\}$ for $j < \ell$.
\end{lemma}

Via the three structural lemmas, the main enumerative result is proved in Theorem~\ref{th:strongn} using a recurrence relation on Stirling numbers of the second kind. 

\begin{theorem}
\label{th:strongn}
There exist $\varphi(\ell,n)$ strongly--connected word--graph members of the family $\wordgraphfam{\ell}{n}$ where:
\begin{align*}
\varphi(\ell,n) &= n!T(\ell, n) \\
T(\ell, n) &= \stirling{\ell-1}{n} + \sum_{j=0}^{\ell-2} \sum_{m=0}^{n-2} \stirling{j}{m} T(\ell-j-1, n-m) (n-m-1) 
\end{align*}
\end{theorem}

The recurrence in Theorem~\ref{th:strongn} was originally derived by Zabrocki\cite{rec} for an algebraic problem but, as yet, no formal proof has been published.

\section{Proofs of results}
In this section we use the previous results from Subsection~\ref{subsec:previous} to prove the theorems and lemmas introduced in Section~\ref{sec:results}. We start with a characterisation of the structure of strongly--connected word--graphs using edge--connectivity.

\label{sec:proofs}
\begin{lemma-edgecon}
Iff $\lambda(G_\omega) > 1$ then $G_\omega$ is strongly--connected.
\end{lemma-edgecon}
\begin{proof}
A strongly--connected digraph $D$, by Definition~\ref{def:strong}, has two edge--disjoint paths between every vertex pair thus, by Menger's Theorem\citep{menger}, $D$ is always 2--edge--connected. The remainder of this proof is dedicated to showing that if $\lambda(G_\omega)> 1$ then $G_\omega$ is strongly--connected.

Let $K_n$ be the complete digraph of order $n$ such that its vertices are labelled with the first $n$ integers. Let $\omega$ be a walk on $K_n$ and let $G_\omega$ be the word--graph represented by $\omega$. We take bidirectional edges in $G_\omega$ to be pairs of unidirectional edges with opposing orientations. If the walk $\omega$ resides in vertex $v_j \in V(K_n)$ then the walk can take one of two actions:
\begin{enumerate}
\renewcommand{\theenumi}{\roman{enumi}}
\item transition from $v_j$ to $v_{j+k}$ for $j+k \le n$ 
\vspace{-0.5em}
\item transition from $v_j$ to $v_{j-k}$ for $j-k \ge 1$ 
\end{enumerate}

Assume, without loss of generality, that $\omega$ starts on $v_1$ and, for the first $n-1$ steps, uses action~i so that $G_\omega$ is weakly--connected (Corollary~\ref{col:weak}) and, by Menger's Theorem\cite{menger}, is 1--edge--connected (Figure~\ref{fig:edge1}). 
\newpage
\begin{figure}[h]
\centering
\begin{tikzpicture}
\centering
\node[state] (A) {\Large $1$};
\node[state,right of=A,draw=white] (B) {\Large $\cdots$};
\node[state,right of=B] (C) {\Large $n$};
\path[->]
(A) edge node {$ $} (B)
(B) edge node {$ $} (C);
\end{tikzpicture}
\caption{$G_\omega$ after $n-1$ steps on $K_n$}
\label{fig:edge1}
\end{figure}
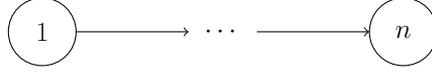

To ensure $G_\omega$ is  2--edge--connected, the walk $\omega$ must continue. The walk cannot use action~i because it already resides in $v_n$ therefore it must use action~ii to make the transition from $v_n$ to $v_{n-j}$ (Figure~\ref{fig:edge2}). 

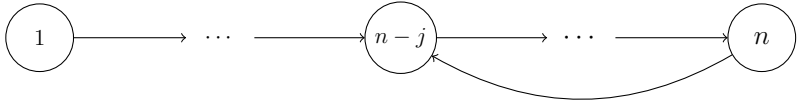
\begin{figure}[h]
\centering
\begin{tikzpicture}
\centering
\node[state] (X) {\large $1$};
\node[state,right of=X,draw=white] (Y) {\large $\cdots$};
\node[state,right of=Y] (A) {\large $n-j$};
\node[state,right of=A,draw=white] (B) {\Large $\cdots$};
\node[state,right of=B] (C) {\Large $n$};
\path[->]
(X) edge node {$ $} (Y)
(Y) edge node {$ $} (A)
(A) edge node {$ $} (B)

(B) edge node {$ $} (C)
(C) edge [bend left] node {$ $} (A);
\end{tikzpicture}
\caption{$G_\omega$ after $n$ steps on $K_n$}
\label{fig:edge2}
\end{figure}

$G_\omega$ now contains the strongly--connected component $C = \{v_n, \ldots, v_{n-j}\}$. By Menger's Theorem\cite{menger}, $C$ is 2--edge--connected. To ensure the entirety of $G_\omega$ is 2--edge--connected, $\omega$ must again make the transition outside $C$ so that $C$ becomes larger still. If $C$ does not encompass all nodes in $K_n$ then $G_\omega$ cannot be 2--edge--connected because $G_\omega$ contains a bridge. Therefore the walk $\omega$ must continue in this manner until the strong--component $C$ encompassed all vertices so that $G_\omega$ is 2--edge--connected and strongly--connected (Figure~\ref{fig:edge3}).

\begin{figure}[ht]
\centering
\begin{tikzpicture}
\centering
\node[state] (X) {\large $1$};
\node[state,right of=X, draw=white] (Y) {\large $\cdots$};
\node[state,right of=Y] (A) {\large $n-j$};
\node[state,right of=A,draw=white] (B) {\Large $\cdots$};
\node[state,right of=B] (C) {\Large $n$};
\path[->]
(X) edge node {$ $} (Y)
(Y) edge node {$ $} (A)
(A) edge node {$ $} (B)

(B) edge node {$ $} (C)
(C) edge [bend left] node {$ $} (A)
(B) edge [bend left] node {$ $} (Y)
(A) edge [bend left] node {$ $} (X);
\end{tikzpicture}
\caption{$G_\omega$ after $n+k$ steps on $K_n$}
\label{fig:edge3}
\end{figure}
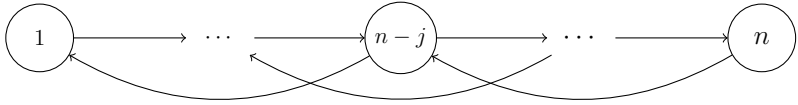

Thus we have shown that 2--edge--connectivity is a necessary and sufficient condition for the strong--connectivity of word--graphs.
\end{proof}
 
Next we show that the structure of words that represent strongly--connected digraphs is characterised by the digraph's \emph{bridgelessness} (2--edge--connectedness).

\begin{lemma-partition}
Iff $\lambda(G_\omega) > 1$ then $\omega$ cannot be partitioned into disjoint factors.
\end{lemma-partition}
\begin{proof}
Let $p$ be a $k-$partition of the $\ell-$word $\omega$. If $p$ is a disjoint partition then $2 \le k \le n$. If $k > n$ then, by the pigeon-hole principle, there must be some pair of non--disjoint partitions because $|A(\omega)| = n$. If $k=1$ then $p_1 = \omega$ and thus cannot be disjoint.

Assume $p$ is a disjoint $k-$partition of $\omega$ for $2 \le k \le n$. For any $p_{i}, p_{i+1}$ partition pair, contiguous in $\omega$, of length $u$ and $v$ respectively, there exists an edge between the vertices $p_{i,u}$ and $p_{i+1,v}$ in $G_\omega$. Given that all $p_i$ partitions are disjoint, vertex $p_{i,u}$ occurs in no other partition but $p_i$ and vertex $p_{i+1,v}$ occurs in no other partition but $p_{i+1}$. Therefore there exists only a single unidirectional edge between $p_{i,u}$ and $p_{i+1,v}$ (Figure~\ref{fig:partition}). 

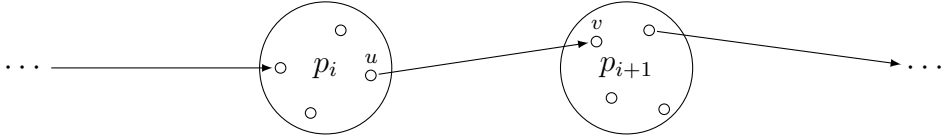
\begin{figure}[!ht]
\centering
\begin{tikzpicture}
\draw (0,0) circle (25pt) node[circle] {$p_i$};
\draw (4,0) circle (25pt) node[circle] {$p_{i+1}$};
\draw[draw=white] (-4,0) circle (25pt) node[circle,draw=white] {$\cdots$};
\draw[draw=white] (8,0) circle (25pt) node[circle,draw=white] {$\cdots$};

\draw (0.2,0.5) circle (2pt) node[circle,] {$ $};
\draw (-0.2,-0.6) circle (2pt) node[circle] {$ $};
\draw (-0.6,0) circle (2pt) node[circle] {$ $};
\draw (4.3,0.5) circle (2pt) node[circle] {$ $};
\draw (4.5,-0.55) circle (2pt) node[circle] {$ $};
\draw (3.8,-0.4) circle (2pt) node[circle] {$ $};

\draw (0.55,0.05) circle (0pt) node {$\ ^u$};
\draw (0.6,-0.1) circle (2pt) node[circle] {$ $};
\draw (3.55,0.50) circle (0pt) node {$\ ^v$};
\draw (3.6,0.35) circle (2pt) node[circle] {$ $};

\begin{scope}[>=latex]
\draw[->,shorten >=3pt, shorten <=3pt] (0.6,-0.1) -- (3.6,0.35);
\draw[->,shorten >=3pt, shorten <=10pt] (-4,0) -- (-0.6,0);
\draw[->,shorten >=10pt, shorten <=3pt] (4.3,0.5) -- (8,0);
\end{scope}
\end{tikzpicture}
\caption{A word--graph represented by a contiguous disjoint partition}
\label{fig:partition}
\end{figure}

Furthermore, because every $p_{i-j}$,$p_{i+1+m}$ partition pair is disjoint, there can be no path between $p_{i,u}$ and $p_{i+1,v}$ except the $\{p_{i,u},p_{i+1,v}\}$ edge. If that edge were to be removed then $G_\omega$ would become disconnected. As such, if $\omega$ can be partitioned into $k > 1$ disjoint factors then $G_\omega$ contains a bridge and $\lambda(G_\omega) = 1$. If $\omega$ cannot be partitioned into disjoint factors then $G_\omega$ cannot contain a bridge and $\lambda(G_\omega)>1$.
\end{proof}

\begin{corollary}
If a word $\omega$ can be minimally partitioned into $k$ disjoint factors then $G_\omega$ contains $k-1$ bridges. 
\end{corollary}

Lemmas~\ref{lem:edgecon} and~\ref{lem:partition} show that words that represent strongly--connected digraphs cannot be partitioned into disjoint factors. Lemma~\ref{lem:subset} restates these results in a manner that facilitates the enumeration of strongly--connected word--graphs.

In Theorem~\ref{th:card}, the $n!$ term of $n! \stirling{\ell}{n}$ accounts for factor permutation. In Lemma~\ref{lem:subset} the $n!$ term will be temporary ignored to make the proof simpler such that we only deal with the set of $\ell-$words over an $n-$alphabet with a unique factorisation, \emph{i.e.} we are not concerned with the order of those factors.

\begin{lemma-subset}
The number of $\ell-$words over an $n-$alphabet, with a unique factorisation, that represent word--graphs with $\lambda(G_\omega) > 1$ is the number of $n-$partitions of an $\ell-$set such that no proper subset of parts has a union equal to the subset $\{1,2,\ldots,j\}$ for $j < \ell$.
\end{lemma-subset}

\begin{proof}
Let $p$ be an $n-$partition of an $\ell-$set with members equal to the first $\ell$ non--zero integers. Let each $p_i$ be mapped to an alphabet member $a_i \in A(\omega)$ of some word $\omega$. Let the $j-$th integer member of $p_{i}$ be the index of letter $a_i \in A(\omega)$ in $\omega$. This construction gives an $\ell-$word over an $n-$alphabet:
\begin{equation}
\label{eq:subsetmap}
\omega_{p_{ij}} = h(p_i) = a_i
\end{equation}

The number of $n-$partitions of an $\ell$-set is $\stirling{\ell}{n}$, which is the same as the number of $\ell-$words over an $n-$alphabet with a unique factorisation.

If a proper $0<k<\ell$ subset of $p$ has a union equal to the subset $\{1, 2, \ldots,j\}$ for $j <\ell$ then $p$ is deemed reducible. A word constructed from a reducible partition is called partition--reducible. For a partition--reducible $\ell-$word, the factor $\omega_1\ldots \omega_j$ is over the alphabet $\{a_1,\ldots,a_k\}$ and the factor $\omega_{j+1}\ldots \omega_\ell$ is over the alphabet $\{a_{k+1},\ldots,a_n\}$. Therefore there exists a disjoint partition of $\omega$ between $\omega_j$ and $\omega_{j+1}$. Thus partition--reducible words cannot represent strongly--connected word--graphs because the word--graphs are \mbox{1--edge--connected} (Lemmas~\ref{lem:edgecon} and~\ref{lem:partition}).

If a partition is not reducible, it is irreducible. Let $\omega$ be a partition--irreducible word formed from the partition $p$ as described by Equation~\ref{eq:subsetmap}. Assume $1 \in p_1$ and because $p$ is irreducible, the union of $p_1$ with any proper subset of $p - \{p_1\}$ cannot be equal to $\{1,\ldots,j\}$. Therefore there exists no disjoint partition of $\omega$. Thus partition--irreducible words always represent strongly--connected word--graphs because the word--graphs are at least 2--edge--connected (Lemmas~\ref{lem:edgecon} and~\ref{lem:partition}).
\end{proof}

Via Lemma~\ref{lem:subset}, a restatement of Lemma~\ref{lem:partition}, we can count the number of strongly--connected word--graphs in $\wordgraphfam{\ell}{n}$.

\begin{theorem-strongn}
There exist $\varphi(\ell,n)$ strongly--connected word--graph members of the family $\wordgraphfam{\ell}{n}$ where:
\begin{align*}
\varphi(\ell,n) &= n!T(\ell, n) \\
T(\ell, n) &= \stirling{\ell-1}{n} + \sum_{j=0}^{\ell-2} \sum_{m=0}^{n-2} \stirling{j}{m}  T(\ell-j-1, n-m) (n-m-1)
\end{align*}
\end{theorem-strongn}
\begin{proof}
There are two ways to construct an irreducible $\ell-$set $n-$partition. First, let $p$ be a partition of  $\{1,2,\ldots,\ell-1\}$ into $n$ subsets. Assume that $1 \in p_1$ and add $\ell$ to $p_1$. If $p$ were a reducible partition then the union of, at maximum, $n-1$ subsets of $p$ would form the set $\{1,\ldots,\ell\}$. Besides $1$ and $\ell$ there exist $\ell-2$ members of the $\ell-$set distributed amongst the subsets of $p$. No subset can be empty thus in order to form the union $\{1,\ldots,\ell\}$, $p_1$ and all remaining $\ell-2$ elements must be in the union. Such a partition cannot be reducible because the union is not formed from a proper subset of the partition. Thus partitions constructed in such a manner are irreducible. There exist:
\begin{align}
\label{eq:partt1}
\stirling{\ell-1}{n}
\end{align}
partitions of this type.

For the second type, let $p$ be an irreducible partition of $\ell-j-1$ into $n-m$ and let $\hat{p}$ be a partition of $j$ into $m$. Let $p\hat{p}$ be the union of the two partitions. The previous irreducible partition type counted all partitions with $\ell \in p_1$ such that $1 \in p_1$. Putting $\ell \in \hat{p}_i$ makes $p\hat{p}$ reducible because $\bigcup p_i = \{1,2,\ldots,\ell-j-1\}$. The only remaining option is to put $\ell$ into $p_i \in p$ such that $1 \notin p_1$, which ensures $p\hat{p}$ is irreducible because: 
\begin{align}
\bigcup& p_i = \{1,2,\ldots,\ell-j-1,\ell\} \notag \\ 
\bigcup& \hat{p}_i = \left(\bigcup p_i\right)^{\complement} = \{\ell-j,\ldots,\ell-1\} \notag \\
\label{eq:notproper}
\bigcup& p\hat{p}_i = \{1,\ldots,\ell\}
\end{align}

Such a partition is irreducible because the only possible union equal to $\{1,\ldots,j\}$ is not formed from a proper subset of $p\hat{p}$ (Equation~\ref{eq:notproper}).

For some fixed $j$ and $m$, there exist $T(\ell-j-1, n-m)$ irreducible $\{1,\ldots,\ell-j-1\}$ partitions into $n-m$; $\stirling{j}{m}$ partitions of $\{\ell-j, \ldots, \ell-1\}$ into $m$; and $n-m-1$ ways of placing element $\ell$ into $p_i \in p$ such that $p$ remains irreducible:
\begin{align}
\label{eq:partt2}
\stirling{j}{m} T(\ell-j-1, n-m) (n-m-1)
\end{align}

To count all such partitions, equation~\ref{eq:partt2} has to be summed over all $0 \le m \le n-2$ because $n-m-1 < 1$ for $m > n-2$ and over all $0 \le j \le \ell-2$ because $\ell-j-1 < 1$ for $j > \ell-2$. 

These two types of irreducible partition constructions must be disjoint because, for the first type, $\ell$ and $1$ are always in the same subset whilst, for the second type, $\ell$ and $1$ are never in the same subset. Likewise, all irreducible partitions are of one of these two types. If $\ell$ and $1$ are in the same subset of $p$ then $p$ is irreducible and must be of the first type because Equation~\ref{eq:partt1} counts all such partitions. If $1 \in p_1$ and $\ell \in p_i$ for $i>1$ and $p$ is irreducible then $p$ is a partition of type two. If $p$ were not of type two then $\ell$ must be in its own subset or in a subset with members of $\{j+1,\cdots,\ell-1\}$ thus the union of a subset of $\{p_1,\ldots,p_{i-1}\}$ would be equal to $\{1,\ldots,j\}$ and so $p$ would be reducible.

Combining Equations~\ref{eq:partt1} and~\ref{eq:partt2} gives the required recurrence relation.
$$T(\ell, n) = \stirling{\ell-1}{n} + \sum_{j=0}^{\ell-2} \sum_{m=0}^{n-2} \stirling{j}{m} T(\ell-j-1, n-m) (n-m-1)$$

The base cases for $T(\ell,n)$ are: 
\begin{align}
\label{eq:base1}
T(\ell,n) = 0 &\text{ for } \ell \le 0 \text{ and } \ell \le n \\
\label{eq:base2}
T(\ell,1) = 1 
\end{align}

The base cases in Equation~\ref{eq:base1} come from Definition~\ref{def:edges}. The base case in Equation~\ref{eq:base2} comes from the fact that the trivial $\ell-$word $\omega$ with $|A(\omega)|=1$ represents the singleton digraph which is strongly--connected. Including the $n!$ term, to account for factor permutation (see Theorem~\ref{th:card}), gives the function $\varphi(\ell, n) = n! T(\ell, n)$ which is the number of strongly--connected word--graphs in $\wordgraphfam{\ell}{n}$.
\end{proof}

\section{Future work}
We have shown how many members of a word--graph family, defined by a particular edge--set encoding with links to stochastic processes, are strongly--connected. Similar reasoning could be used for the number of word--graph with $k$ strong--components.

\begin{conjecture}
The number of $\ell-$words over an $n-$alphabet, with a unique factorisation, that represent word--graphs with $k$ strong--components is the number of $n-$partitions of an $\ell-$set with $k$ disjoint proper subsets of parts which have a union equal to the subset $\{j,j+1,\ldots,j+m\}$. Each such union is an alphabet--disjoint partition in the corresponding representational word and is a strong--component in the corresponding word--graph.
\end{conjecture}

The word--graph family $\wordgraphfam{\ell}{n}$ includes isomorphic digraphs. The set of $\ell-$words over an $n-$alphabet that represent non--isomorphic words--graphs can be thought of as a \emph{fundamental} word--graph family. A standing question is: what is the cardinality of an arbitrary fundamental word--graph family and for what $\ell$, when $n$ is fixed, is the cardinality maximised? These topics are related to the concept of minimal representational word length\cite{wordgraph}.

\section*{Acknowledgements}
The authors wish to thank Mike Zabrocki\footnote{zabrocki@mathstat.yorku.ca} of York University Canada for providing the recurrence relation in Theorem~\ref{th:strongn} and for giving a sketch of the proof, and thank the ESRC for funding this research.

\bibliographystyle{plainnat}
\bibliography{strong-refs}

\begin{thebibliography}{8}
\providecommand{\natexlab}[1]{#1}
\providecommand{\url}[1]{\texttt{#1}}
\expandafter\ifx\csname urlstyle\endcsname\relax
  \providecommand{\doi}[1]{doi: #1}\else
  \providecommand{\doi}{doi: \begingroup \urlstyle{rm}\Url}\fi

\bibitem[Alspach and Pullman(1974)]{path74}
B.~Alspach and N.~J. Pullman.
\newblock Path decomposition of digraphs.
\newblock \emph{Bulletin of the Australian Mathematical Society.}, 10:\penalty0
  421--427, 1974.

\bibitem[Bell et~al.(In Review{\natexlab{a}})Bell, Rayson, and
  Berridge]{word-graph-grammar}
Edward J.~L. Bell, Paul Rayson, and Damon Berridge.
\newblock Word--graph enumeration with context--free grammars.
\newblock In Review{\natexlab{a}}.

\bibitem[Bell et~al.(In Review{\natexlab{b}})Bell, Rayson, and
  Berridge]{wordgraph}
Edward J.~L. Bell, Paul Rayson, and Damon Berridge.
\newblock On the structural characterisation of word--graphs.
\newblock In Review{\natexlab{b}}.

\bibitem[Kitaev and Pyatkin(2008)]{representable}
Sergey Kitaev and Artem Pyatkin.
\newblock On representable graphs.
\newblock \emph{Journal of Automata, Languages and Combinatorics.}, 13\penalty0
  (1):\penalty0 45--54, 2008.

\bibitem[Li and Zhang(2010)]{overlap}
Xianyue Li and Heping Zhang.
\newblock Embedding on alphabet overlap digraphs.
\newblock \emph{Journal of Mathematical Chemistry.}, 47\penalty0 (1):\penalty0
  62--71, 2010.

\bibitem[Menger(1927)]{menger}
Karl Menger.
\newblock Zur allgemeinen kurventheorie.
\newblock \emph{Fundamenta Mathematicae}, 10:\penalty0 95--115, 1927.

\bibitem[Shur(2010)]{word-tern}
Arseny~M Shur.
\newblock On ternary square-free circular words.
\newblock \emph{The Electronic Journal of Combinatorics.}, 17\penalty0 (1),
  2010.

\bibitem[Sloane(2003)]{rec}
N.~J.~A Sloane, 2003.
\newblock Sequence in A087903 \emph{The On-Line Encyclopedia of Integer
  Sequences} -- \url{http://oeis.org/A087903}. Submitted by Mike Zabrocki.

\end{thebibliography}
\end{document}